# Existence and Stability of Solutions to Non-Lipschitz Stochastic Differential Equations Driven by Lévy Noise

Yong Xu[1)], Bin Pei

Department of Applied Mathematics, Northwestern Polytechnical University, Xi'an, 710072, China

**Abstract:** In this paper, the successive approximation method is applied to investigate the existence and uniqueness of solutions to the stochastic differential equations (SDEs) driven by Lévy noise under non-Lipschitz condition which is a much weaker condition than Lipschiz one. The stability of the solutions to non-Lipschitz SDEs driven by Lévy noise is also considered, and the stochastic stability is obtained in the sense of mean square.

## 1. Introduction

In the real world, random fluctuations appear commonly in various natural and synthetic systems, and stochastic differential equations (SDEs) with random fluctuations and noises have been applied as the mathematical models of many practical systems. Such models can describe a great deal of different scientific and engineering problems, which appear in different fields including biology, epidemiology, mechanic, economics, finance and so forth [1-5]. So it is natural and necessary to investigate dynamical properties of the solutions to SDEs to find the effects of random perturbations in the corresponding realistic systems.

The mathematical models obtained have been greatly developed for SDEs under a random disturbance of the "Gaussian white noise", namely, the investigations concerning SDEs driven by Brownian motion have been very abundant up to now [4, 27]. However, Gaussian noises just like Brownian motion are not always appropriate while there exist large external and/or internal fluctuations with possible big jumps in some practical situations and environments. Unlike Brownian motion with the properties of no jump (the sample paths are continuous almost surely), normal diffusion (the mean square displacement increases linearly in time), and light tail or exponential relaxation (the probability density function decays exponentially in space) [4], non-Gaussian Lévy noise has completely peculiar properties of anomalous diffusion (mean square displacement is a nonlinear power law of time) [6] and heavy tail distribution or non-exponential relaxation [7]. Given the above-mentioned peculiar properties, non-Gaussian Lévy noise may be more appropriate to model the stochastic disturbances especially with extremely large jumps [8-11, 28-30]. For instance, recently, there has a growing interest in the use of Lévy process to model market behaviors in finance; not only are these of great mathematical interest but also there is growing evidence that they may be more realistic models than those that insist on continuous samples paths [12-14]. Therefore, it is significant to investigate the SDEs, properties of solutions and the applications with non-Gaussian Lévy noise.

---

[1] Corresponding author. E-mail:hsux3@nwpu.edu.cn.
Tel./Fax. : 86-29 8843 1637

A number of meritorious results concerning SDEs with non-Gaussian Lévy noise has been presented in existing literatures [7, 8, 15-17]. Among them, conditions which can guarantee the existence and uniqueness of the solutions to the SDEs with non-Gaussian Lévy noise are to be assumed as the one of the most basic and important Lipschitz condition. Generally, the Lipschitz case is a very common condition which has been widely used, and the existence and uniqueness of solutions to SDEs with non-Gaussian Lévy noise has been proved [16]. But, this condition is seemed to be considerably trenchant when one discusses variable applications in real world [7]. Thus, the importance to find some weaker and more generalized condition under which the SDEs with non-Gaussian Lévy noise still exist unique solutions is self-evident.

Fortunately, the so-called non-Lipschitz conditions have been proposed in [18-23], and by the method of the successive approximation, we in this paper prove the existence and uniqueness of solutions to SDEs driven by Lévy noise under this kind of non-Lipschitz condition with Lipschitz one as a special case. Furthermore, in order to obtain the dynamic properties of the solutions to non-Lipschitz SDEs driven by Lévy noise, we present the stability conditions and the theorem that the solutions are stochastic stable in the sense of mean square.

## 2. Some preparations

### 2.1 Lévy motions

**Definition 1.** Let $L(t)$ be a stochastic process defined on a probability space $(\Omega, \mathcal{F}, P)$, if

(i) $L(0) = 0$ (a.s.),

(ii) $L(t)$ has independent and stationary increments,

(iii) $L(t)$ is stochastically continuous, i.e. for all $a > 0$ and for all $s \geq 0$,

$$\lim_{t \to s} P(|L(t) - L(s)| > a) = 0,$$

then $L(t)$ is a Lévy process.

A scalar Lévy process or motion is characterized by a drift parameter $q \in R^d$, a covariance $d \times d$ matrix $A$ and a nonnegative Borel measure $v$ defined on $(R^d, \mathcal{B}(R^d))$ and concentrated on $R^d \setminus \{0\}$, which satisfies

$$\int_{R^d \setminus \{0\}} (x^2 \wedge 1) v(dx) < \infty.$$

This measure $v$ is the so called Lévy jump measure of the Lévy process $L(t)$, and $(A, v, q)$ is defined as the generating triplet of Lévy motion.

It is known to all that a scalar Lévy motion is completely determined by the Lévy-Khintchine formula [8, 16, 24]. Now, we introduce the Lévy-Itô decomposition theorem with no proof and we may refer to [16] for more details.

**Proposition 1.** If $L(t)$ is a Lévy process in $R^d$, there exists $b_1 \in R^d$, a Brownian

motion $B(t)$ and independent Poisson random measure $N$ on $R_+ \times (R^d - \{0\})$, such that for each $t \geq 0$,

$$L(t) = b_1 t + B(t) + \int_{|x|<c} x \tilde{N}(t, dx) + \int_{|x|\geq c} x N(t, dx), \qquad (1)$$

where $N(dt, dx)$ is Poisson random measure, $\tilde{N}(dt, dx) = N(dt, dx) - v(dx) dt$ is the compensated Poisson random measure for $L(t)$ and $b_1 = E\left(L(1) - \int_{|x|\geq c} x N(1, dx)\right)$ with the parameter $c \in [0, \infty)$ a constant.

## 2.2 SDEs driven by Lévy noise

We concern the SDEs with Lévy noise on $R^d$:

$$dX(t) = f(t, X(t-)) dt + g(t, X(t-)) dL(t), t \in [0, T]. \qquad (2)$$

Using the Lévy-Itô decomposition (1), we can rewrite this as follows

$$dX(t) = f_1(t, X(t-)) dt + g(t, X(t-)) dB(t) + \int_{|x|<c} g(t, X(t-)) x \tilde{N}(dt, dx) \\ + \int_{|x|\geq c} g(t, X(t-)) x N(dt, dx), \qquad (3)$$

where $f_1 = b_1 + f$.

Hence Eq. (2) can be rewritten in the following more general form:

$$dX(t) = b(t, X(t-)) dt + \sigma(t, X(t-)) dB(t) + \int_{|x|<c} F(t, X(t-), x) \tilde{N}(dt, dx) \\ + \int_{|x|\geq c} G(t, X(t-), x) N(dt, dx). \qquad (4)$$

The third and fourth terms in the right hand side of Eq. (4) represent small and large jumps, and the term involving $G$ (or $F$ respectively) is absent when we take $c = \infty$. The term in Eq. (4) involving large jumps is controlled by $G$, which can be neglected through the technique of interlacing [16]. So we only focus on the study of an equation driven by continuous noise interspersed with small jumps, and then Eq. (4) can be modified as

$$dX(t) = b(t, X(t-)) dt + \sigma(t, X(t-)) dB(t) + \int_{|x|<c} F(t, X(t-), x) \tilde{N}(dt, dx). \quad (5)$$

Then the integral form of Eq. (5) can be given

$$X(t) = X(0) + \int_0^t b(s, X(s-)) ds + \int_0^t \sigma(s, X(s-)) dB(s)$$
$$+ \int_0^t \int_{|x|<c} F(s, X(s-), x) \tilde{N}(ds, dx), \tag{6}$$

with initial value $X(0) = \xi$, where $\xi$ is a given $d$-dimensional random vector. $b(t,\bullet)$ and $F(t,\bullet,x)$ are given $d$-dimensional random vector functions. $\sigma(t,\bullet)$ is a $d \times r$ matrix, $c$ is a positive constant, and $t \in [0,T]$, $0 < T < +\infty$. $B(t)$ represents $r$-dimensional Brownian motion. We define the components of these terms as follows. For each $i$ and $j$, the mappings $b^i : [0,T] \times R^d \to R$, $\sigma_j^i : [0,T] \times R^d \to R$ and $F^i : [0,T] \times R^d \times R^d \to R$ are all assumed to be measurable for $1 \leq i \leq d, 1 \leq j \leq r$.

## 2.3 Non-Lipschitz condition

Now, consider the following assumption on the coefficients of Eq. (6):

**Assumption1.** For each fixed $t \in [0,T]$, let $b(t,y)$, $\sigma(t,y)$ and $F(t,y,x)$ be continuous in $y$, and $\forall y_1, y_2 \in R^n$,

$$\begin{aligned}&|b(t,y_1) - b(t,y_2)|^2 + \|\sigma(t,y_1) - \sigma(t,y_2)\|^2 \\ &+ \int_{|x|<c} |F(t,y_1,x) - F(t,y_2,x)|^2 v(dx) \leq \lambda(t) \kappa\left(|y_1 - y_2|^2\right),\end{aligned} \tag{7}$$

where

(1a) $\lambda(t) : [0,\infty) \to R^+$ is an integrable function, $\kappa(q)$ is monotone non-decreasing, continuous and $\kappa(q)$ or $\kappa(q)^2/q$ is a concave function with respect to $q$ for fixed $t \geq 0$ with $\kappa(0) = 0$ such that

$$\int_{0+} \frac{1}{\kappa(q)} dq = \infty,$$

(1b) $$\|\sigma(t,y_1) - \sigma(t,y_2)\|^2 = \sum_{i=1}^d \sum_{j=1}^r \left(\sigma_j^i(t,y_1) - \sigma_j^i(t,y_2)\right)^2,$$

$$|b(t,y_1) - b(t,y_2)|^2 = \sum_{i=1}^d \left(b^i(t,y_1) - b^i(t,y_2)\right)^2,$$

$$|F(t,y_1,x)-F(t,y_2,x)|^2 = \sum_{i=1}^{d}(F^i(t,y_1,x)-F^i(t,y_2,x))^2,$$

(1c) $\sigma(t,0), b(t,0), F(t,0,x)$ are integrable.

The Assumption 1 is the so-called non-Lipschitz conditions. In particular, if $\kappa(q)=q, \lambda(t)=K$, then the Assumption 1 reduces to Lipschitz condition. In other words, non-Lipschitz condition is weaker than Lipschitz condition.

**Remark 1.** For fixed $t \geq 0$, if $\kappa(q)$ or $\kappa(q)^2/q$ is concave function with respect to $q$, and $\forall y_1, y_2 \in R^n$, the following inequality

$$E|b(t,y_1)-b(t,y_2)|^2 + E\|\sigma(t,y_1)-\sigma(t,y_2)\|^2 \\ + E\int_{|x|<c}|F(t,y_1,x)-F(t,y_2,x)|^2 v(dx) \leq \lambda(t)\kappa(E|y_1-y_2|^2), \quad (8)$$

holds.

The proof is given in [18, 23]. □

## 3. Existence and uniqueness of solutions to non-Lipschitz SDEs with Lévy noise

We in this section consider the existence and uniqueness of solutions to SDEs driven by Lévy noise under the above-mentioned non-Lipschitz condition. Hereafter, we assume without losing generality that $1 \leq T < +\infty$, and $K_1, C_i, i=1,2\cdots 6$ are all constants.

Define a sequence of stochastic process $\{X_k(t)\}_{k=1,2\cdots}$ with $X_0(t) \equiv \xi$ a random variable and $E|\xi|^2 < +\infty$ as follows:

$$X_k(t) = \xi + \int_0^t b(s, X_{k-1}(s-))ds + \int_0^t \sigma(s, X_{k-1}(s-))dB(s) \\ + \int_0^t \int_{|x|<c} F(s, X_{k-1}(s-), x)\tilde{N}(ds,dx), k=1,2\cdots \quad (9)$$

### 3.1 Existence of solutions to SDEs

**Theorem 1.** Suppose that $b(t,y)$, $F(t,y,x)$ and $\sigma(t,y)$ satisfy Assumption 1, then

$$\lim_{n,i\to\infty} E\left(\sup_{0\le t\le T} |X_n(t) - X_i(t)|^2\right) = 0. \tag{10}$$

By Theorem 1, we say $\{X_k(\bullet)\}_{k=1,2\cdots}$ is a Cauchy sequence with its limit $X(\bullet)$. Then letting $n\to\infty$ in (9), we finally obtain that the solutions to Eq. (2) exist.

To prove the Theorem 1, several auxiliary assertions are required.

**Lemma 1. (Doob's martingale inequality)** If $(X(t), t\ge 0)$ is a positive submartingale, then for any $p>1, t>0$,

$$E\left[\sup_{0\le s\le t} X(s)^p\right] \le \left(\frac{p}{p-1}\right)^p E(X(t)^p).$$

Refer to [16] to obtain the proof. □

**Lemma 2.** Under Assumption 1, there exists a positive number $K_1 > 0$ such that $\forall (t,y)\in[0,T]\times R^d$,

$$|b(t,y)|^2 + \|\sigma(t,y)\|^2 + \int_{|x|<c} |F(t,y,x)|^2 v(dx) \le K_1(1+|y|^2). \tag{11}$$

**Proof:** Since $\kappa(q)$ or $\kappa(q)^2/q$ is a concave and non-negative function, such that $\kappa(0)=0$, we can choose two positive constants $a>0$ and $b>0$, so that [18]

$$\kappa(q) \le a + bq, q\ge 0.$$

Then, using (7) to yields

$$|b(t,y)|^2 + \|\sigma(t,y)\|^2 + \int_{|x|<c} |F(t,y,x)|^2 v(dx)$$
$$\le 2|b(t,y)-b(t,0)|^2 + 2|b(t,0)|^2 + 2\|\sigma(t,y)-\sigma(t,0)\|^2 + 2\|\sigma(t,0)\|^2$$
$$+ 2\int_{|x|<c} |F(t,y,x)-F(t,0,x)|^2 vdx + 2\int_{|x|<c} |F(t,0,x)|^2 v(dx)$$
$$\le 2\left(|b(t,y)-b(t,0)|^2 + \|\sigma(t,y)-\sigma(t,0)\|^2 + \int_{|x|<c} |F(t,y,x)-F(t,0,x)|^2 v(dx)\right)$$

$$+2|b(t,0)|^2 + 2\|\sigma(t,0)\|^2 + 2\int_{|x|<c}|F(t,0,x)|^2 v(dx)$$

$$\leq 2\sup_{0\leq t\leq T}\left(|b(t,0)|^2 + \|\sigma(t,0)\|^2 + \int_{|x|<c}|F(t,0,x)|^2 v(dx)\right) + 2\lambda(t)\kappa(|y|^2)$$

$$\leq K_1(1+|y|^2),$$

where

$$K_1 = \max_{0\leq t\leq T}\left\{2\sup_{0\leq t\leq T}\left(|b(t,0)|^2 + \|\sigma(t,0)\|^2 + \int_{|x|<c}|F(t,0,x)|^2 v(dx) + \lambda(t)a\right), 2b\lambda(t)\right\} < +\infty.$$

This completes the proof of Lemma 2. □

**Lemma 3.** Under Lemma 2, one can have

$$E|X_k(t)|^2 \leq C_1, t \in [0,T], k = 1,2\cdots \qquad (12)$$

**Proof :** Here, for $k = 1, 2\cdots$, we shall show

$$E|X_k(t)|^2 \leq 4E|\xi|^2 \sum_{l=0}^{k}\frac{(4K_1T)^l}{l!}t^l + \sum_{l=1}^{k}\frac{(4K_1T)^l}{l!}t^l. \qquad (13)$$

Firstly, for $k = 1$, it is easy to verify that

$$E|X_1(t)|^2$$

$$= E\left|\xi + \int_0^t b(s,X_0(s-))ds + \int_0^t \sigma(s,X_0(s-))dB(s) + \int_0^t\int_{|x|<c}F(s,X_0(s-),x)\tilde{N}(ds,dx)\right|^2$$

$$\leq 4E|\xi|^2 + 4TE\int_0^t\left(|b(s,X_0(s-))|^2 + \|\sigma(s,X_0(s-))\|^2 + \int_{|x|<c}|F(s,X_0(s-),x)|^2 v(dx)\right)ds$$

$$\leq 4E|\xi|^2 + 4K_1T\int_0^t(1+E|X_0(s-)|^2)ds$$

$$\leq 4E|\xi|^2 + 4K_1T(1+E|\xi|^2)t.$$

So, (13) holds for $k = 1$.

Secondly, assume that (13) holds for $k$, then, by induction, we have for $k+1$,

$$E|X_{k+1}(t)|^2$$

$$= E\left|\xi + \int_0^t b(s,X_k(s-))ds + \int_0^t \sigma(s,X_k(s-))dB(s) + \int_0^t\int_{|x|<c}F(s,X_k(s-),x)\tilde{N}(ds,dx)\right|^2$$

$$\leq 4E|\xi|^2 + 4TE\int_0^t\left(|b(s,X_k(s-))|^2 + \|\sigma(s,X_k(s-))\|^2 + \int_{|x|<c}|F(s,X_k(s-),x)|^2 v(dx)\right)ds$$

$$\leq 4E|\xi|^2 + 4K_1 T \int_0^t \left(1 + 4E|\xi|^2 \sum_{l=0}^{k} \frac{(4K_1 T)^l}{l!} t^l + \sum_{l=1}^{k} \frac{(4K_1 T)^l}{l!} t^l \right) ds$$

$$= 4E|\xi|^2 + 4K_1 Tt + 4E|\xi|^2 \sum_{l=1}^{k+1} \frac{(4K_1 T)^l}{l!} t^l + \sum_{l=2}^{k+1} \frac{(4K_1 T)^l}{l!} t^l$$

$$= 4E|\xi|^2 \sum_{l=0}^{k+1} \frac{(4K_1 T)^l}{l!} t^l + \sum_{l=1}^{k+1} \frac{(4K_1 T)^l}{l!} t^l.$$

Therefore, (13) holds for all $k$. Now, from (13), we obtain

$$E|X_k(t)|^2 \leq 4\left(1 + E|\xi|^2\right)\exp\{4K_1 T^2\}. \tag{14}$$

This completes the proof of Lemma 3. □

**Lemma 4.** Suppose that $b(t,y), F(t,y,x)$ and $\sigma(t,y)$ satisfy the Assumption 1, then for $t \in [0,T]$, $n \geq 1, k \geq 1$,

$$E\left(\sup_{0 \leq s \leq t} |X_{n+k}(s) - X_n(s)|^2\right)$$
$$\leq C_2 \int_0^t \kappa \left[E\left(\sup_{0 \leq s_1 \leq s} |X_{n+k-1}(s_1-) - X_{n-1}(s_1-)|^2\right)\right] ds, \tag{15}$$

and

$$E\left(\sup_{0 \leq s \leq t} |X_{n+k}(s) - X_n(s)|^2\right) \leq C_3 t. \tag{16}$$

**Proof:** By Lemmas 1, 2 and Remark 1, we obtain

$$E\left(\sup_{0 \leq s \leq t} |X_{n+k}(s) - X_n(s)|^2\right)$$

$$\leq 3E\left\{\sup_{0 \leq s \leq t} \left|\int_0^s \left[b(s_1, X_{n+k-1}(s_1-)) - b(s_1, X_{n-1}(s_1-))\right] ds_1\right|^2\right\}$$

$$+ 3E\left\{\sup_{0 \leq s \leq t} \left|\int_0^s \left[\sigma(s_1, X_{n+k-1}(s_1-)) - \sigma(s_1, X_{n-1}(s_1-))\right] dB(s_1)\right|^2\right\}$$

$$+ 3E\left\{\sup_{0 \leq s \leq t} \left|\int_0^s \int_{|x|<c} \left[F(s_1, X_{n+k-1}(s_1-), x) - F(s_1, X_{n-1}(s_1-), x)\right] \tilde{N}(ds_1, dx)\right|^2\right\}$$

$$\leq 3TE\int_0^t \left|b(s, X_{n+k-1}(s-)) - b(s, X_{n-1}(s-))\right|^2 ds$$

$$+12E\int_0^t \left\|\sigma(s, X_{n+k-1}(s-)) - \sigma(s, X_{n-1}(s-))\right\|^2 ds$$

$$+12E\int_0^t \int_{|x|<c} \left|F(s, X_{n+k-1}(s-), x) - F(s, X_{n-1}(s-), x)\right|^2 v(dx)ds$$

$$\leq 12T\left(\sup_{0\leq s\leq T}\lambda(s)\right)\int_0^t \kappa\left[E\left(\sup_{0\leq s_1\leq s}\left|X_{n+k-1}(s_1-) - X_{n-1}(s_1-)\right|^2\right)\right]ds.$$

So, set $C_2 = 12T\left(\sup_{0\leq s\leq T}\lambda(s)\right)$, the inequality (15) holds.

Next, from Lemma 2, 3,

$$E\left(\sup_{0\leq s\leq t}\left|X_{n+k}(s) - X_n(s)\right|^2\right) \leq C_2\int_0^t \kappa\left(E\left(\sup_{0\leq s_1\leq s}\left|X_{n+k-1}(s_1) - X_{n-1}(s_1)\right|^2\right)\right)ds$$

$$\leq C_2\int_0^t \kappa(4K_1)ds.$$

Thus

$$E\left(\sup_{0\leq s\leq t}\left|X_{n+k}(s) - X_n(s)\right|^2\right) \leq C_3 t.$$

This completes the proof of Lemma 4.   □

Now, we choose $0 < T_1 \leq T$, $t \in [0, T_1]$, such that

$$\kappa_1(C_3 t) \leq C_3, \tag{17}$$

where

$$\kappa_1(q) = C_2 \kappa(q).$$

Then, for fixed $k \geq 1$, we introduce two sequences of functions $\{\phi_n(t)\}_{n=1,2...}$ and $\{\tilde{\phi}_{n,k}(t)\}_{n=1,2...}$, where

$$\phi_1(t) = C_3 t,$$

$$\phi_{n+1}(t) = \int_0^t \kappa_1(\phi_n(s))ds,$$

$$\tilde{\phi}_{n,k}(t) = E\left[\sup_{0\leq s\leq t}\left|X_{n+k}(s) - X_n(s)\right|^2\right], n = 1, 2...$$

**Lemma 5.** In terms of Assumption 1,

$$0 \leq \tilde{\phi}_{n,k}(t) \leq \phi_n(t) \leq \phi_{n-1}(t) \leq \cdots \leq \phi_1(t), t \in [0, T_1], \quad (18)$$

holds for each $k \geq 1$ and all positive integer $n \geq 1$.

**Proof:** Firstly, we shall show (18) for $n = 1$. By Lemma 4, we have

$$\tilde{\phi}_{1,k}(t) = E\left(\sup_{0 \leq s \leq t} |X_{1+k}(s) - X_1(s)|^2\right) \leq C_3 t = \phi_1(t), t \in [0, T_1].$$

Secondly, it is easy to verify

$$\tilde{\phi}_{2,k}(t) = E\left(\sup_{0 \leq s \leq t} |X_{2+k}(s) - X_2(s)|^2\right)$$

$$\leq C_2 \int_0^t \kappa\left[E\left(\sup_{0 \leq s_1 \leq s} |X_{1+k}(s_1-) - X_1(s_1-)|^2\right)\right] ds$$

$$\leq \int_0^t \kappa_1\left(\tilde{\phi}_{1,k}(s)\right) ds \leq \int_0^t \kappa_1\left(\phi_1(s)\right) ds = \phi_2(t)$$

$$\leq \int_0^t \kappa_1(C_3 s) ds \leq C_3 t = \phi_1(t),$$

thus, for $n = 2$,

$$\tilde{\phi}_{2,k}(t) \leq \phi_2(t) \leq \phi_1(t), t \in [0, T_1].$$

At last, by what assumed for $n$, for $n+1$, it is easy to verify

$$\tilde{\phi}_{n+1,k}(t) = E\left(\sup_{0 \leq s \leq t} |X_{n+k+1}(s) - X_{n+1}(s)|^2\right)$$

$$\leq C_2 \int_0^t \kappa\left[E\left(\sup_{0 \leq s_1 \leq s} |X_{n+k}(s_1-) - X_n(s_1-)|^2\right)\right] ds$$

$$\leq \int_0^t \kappa_1\left(\tilde{\phi}_{n,k}(s)\right) ds \leq \int_0^t \kappa_1\left(\phi_n(s)\right) ds = \phi_{n+1}(t)$$

$$\leq \int_0^t \kappa_1\left(\phi_{n-1}(s)\right) ds = \phi_n(t).$$

This completes the proof of Lemma 5. □

**The proof of Theorem 1:**

**Step 1:** Since $\phi_n(t)$ decreases monotonically when $n \to \infty$ and $\phi_n(t)$ is non-negative function on $t \in [0, T_1]$, according to Lemma 5, we can define the function $\phi(t)$ by $\phi_n(t) \downarrow \phi(t)$. Obviously, $\phi(0) = 0$ and $\phi(t)$ is continuous function on $[0, T_1]$. Using the definition of $\phi_n(t)$ and $\phi(t)$, we reach

$$\phi(t) = \lim_{k \to \infty} \phi_{n+1}(t) = \lim_{k \to \infty} \int_0^t \kappa_1(\phi_n(s)) ds = \int_0^t \kappa_1(\phi(s)) ds, t \in [0, T_1]. \quad (19)$$

Then, (19) implies $\phi(t) \equiv 0, t \in [0, T_1]$, since $\phi(0) = 0$ and

$$\int_{0+} \frac{dq}{\kappa_1(q)} = \frac{1}{C_2} \int_{0+} \frac{dq}{\kappa(q)} = +\infty.$$

Therefore we obtain

$$0 \leq \lim_{k,n \to \infty} E\left( \sup_{0 \leq t \leq T_1} |X_{n+k}(t) - X_n(t)|^2 \right) = \lim_{k,n \to \infty} \tilde{\phi}_{n,k}(T_1) \leq \lim_{n \to \infty} \phi_n(T_1) = 0,$$

namely,

$$\lim_{n,i \to \infty} E\left( \sup_{0 \leq t \leq T_1} |X_n(t) - X_i(t)|^2 \right) = 0. \quad (20)$$

This completes the Step 1. □

**Step 2**: Define

$$T_2 = \sup\left\{ \tilde{T}: \ \tilde{T} \in [0, T] \text{ and } \lim_{n,i \to \infty} E\left( \sup_{0 \leq t \leq \tilde{T}} |X_n(t) - X_i(t)|^2 \right) = 0 \right\}.$$

Immediately, we get $0 < T_1 \leq T_2 \leq T$. Let $\varepsilon > 0$ be an arbitrary positive number.

Choose $\delta$ so that $0 < \delta < \min(T_2, 1)$ and

$$C_4 \delta < \frac{\varepsilon}{10}. \quad (21)$$

On one hand, from the definition of $T_2$, we have

$$\lim_{n,i \to \infty} E\left( \sup_{0 \leq t \leq T_2 - \delta} |X_n(t) - X_i(t)|^2 \right) = 0.$$

Thus, for large enough $N$, we observe

$$E\left( \sup_{0 \leq t \leq T_2 - \delta} |X_n(t) - X_i(t)|^2 \right) < \frac{\varepsilon}{10}, n, i \geq N. \quad (22)$$

On the other hand, we obtain

$$E\left( \sup_{T_2 - \delta \leq t \leq T_2} |X_n(t) - X_i(t)|^2 \right) \leq 3I_1 + 3I_2 + 3I_3,$$

where

$$I_1 = E\left(\sup_{T_2-\delta \leq t \leq T_2} |X_n(t) - X_n(T_2-\delta)|^2\right),$$

$$I_2 = E\left(|X_n(T_2-\delta) - X_i(T_2-\delta)|^2\right),$$

$$I_3 = E\left(\sup_{T_2-\delta \leq t \leq T_2} |X_i(T_2-\delta) - X_i(t)|^2\right).$$

Using Schwarz's inequality and Lemma1, we get

$$I_1 \leq 3\delta E \int_{T_2-\delta}^{T_2} |b(s_1, X_{n-1}(s_1-))|^2 ds_1 + 12E \int_{T_2-\delta}^{T_2} \|\sigma(s_1, X_{n-1}(s_1-))\|^2 ds_1$$

$$+ 12E \int_{T_2-\delta}^{T_2} \int_{|x|<c} |F(s_1, X_{n-1}(s_1-), x)|^2 v(dx) ds_1$$

$$\leq 12\delta \int_{T_2-\delta}^{T_2} K_1\left(1 + E|X_{n-1}(s_1-)|^2\right) ds_1$$

$$\leq \delta C_4,$$

where $C_4 = 12\delta K_1(1+C_1)$.

Therefore, by (21) we obtain

$$I_1 \leq \frac{\varepsilon}{10},$$

and

$$I_3 \leq \frac{\varepsilon}{10}.$$

In addition, (22) implies

$$I_2 = E|X_k(T_2-\delta) - X(T_2-\delta)|^2 < \frac{\varepsilon}{10}, k \geq N.$$

Now, it is easy to verify

$$E\left(\sup_{0 \leq t \leq T_2} |X_n(t) - X_i(t)|^2\right)$$

$$\leq E\left(\sup_{0 \leq t \leq T_2-\delta} |X_n(t) - X_i(t)|^2\right) + E\left(\sup_{T_2-\delta \leq t \leq T_2} |X_n(t) - X_i(t)|^2\right)$$

$$\leq \frac{\varepsilon}{10} + 3I_1 + 3I_2 + 3I_3 < \varepsilon,$$

namely

$$\lim_{n,i \to \infty} E\left(\sup_{0 \leq t \leq T_2} |X_n(t) - X_i(t)|^2\right) = 0.$$

This completes the Step 2. □

**Step 3:** In this step we shall use the method of reduction to absurdity to show $T_2 = T$.

Firstly, we give two lemmas:

**Lemma 6.** Assume $T_2 < T$ and choose a sequence of number $\{a_i\}_{1 \leq i < n}$ so that $a_i \downarrow 0 (i \to +\infty)$,

$$E\left[\sup_{0 \leq t \leq T_2} |X_n(t) - X_i(t)|^2\right] \leq a_i, \tag{23}$$

then, for $n > i \geq 1$ we have

$$E\left[\sup_{T_2 \leq s \leq T_2 + t} |X_n(s) - X_i(s)|^2\right] \leq 4a_{i-1} + C_5 t, T_2 + t \leq T. \tag{24}$$

**Proof:** First

$$E\left(\sup_{T_2 \leq s \leq T_2 + t} |X_n(s) - X_i(s)|^2\right) \leq 4E|X_{n-1}(T_2) - X_{i-1}(T_2)|^2$$

$$+ 4E\left|\sup_{T_2 \leq s \leq T_2 + t} \int_{T_2}^s [b(s_1, X_{n-1}(s_1-)) - b(s_1, X_{i-1}(s_1-))] ds_1\right|^2$$

$$+ 4E\left|\sup_{T_2 \leq s \leq T_2 + t} \int_{T_2}^s [\sigma(s_1, X_{n-1}(s_1-)) - \sigma(s_1, X_{i-1}(s_1-))] dB(s_1)\right|^2$$

$$+ 4E\left|\sup_{T_2 \leq s \leq T_2 + t} \int_{T_2}^s \int_{|x|<c} [F(s_1, X_{n-1}(s_1-), x) - F(s_1, X_{i-1}(s_1-), x)] \tilde{N}(ds_1, dx)\right|^2$$

$$\leq 4a_{i-1} + 4TE \int_{T_2}^{T_2+t} |b(s_1, X_{n-1}(s_1-)) - b(s_1, X_{i-1}(s_1-))|^2 ds_1$$

$$+ 16E \int_{T_2}^{T_2+t} \|\sigma(s_1, X_{n-1}(s_1-)) - \sigma(s_1, X_{i-1}(s_1-))\|^2 ds_1$$

$$+ 16E \int_{T_2}^{T_2+t} \int_{|x|<c} |F(s_1, X_{n-1}(s_1-), x) - F(s_1, X_{i-1}(s_1-), x)|^2 v(dx) ds_1$$

$$\leq 4a_{i-1} + 16T\left(\sup_{0 \leq s \leq T} \lambda(s)\right) \int_{T_2}^{T_2+t} \kappa\left[E\left(\sup_{T_2 \leq s_1 \leq s} |X_{n-1}(s_1-) - X_{i-1}(s_1-)|^2\right)\right] ds$$

$$\leq 4a_{i-1} + C_6 \int_{T_2}^{T_2+t} \kappa\left[E\left(\sup_{T_2 \leq s_1 \leq s} |X_{n-1}(s_1-) - X_{i-1}(s_1-)|^2\right)\right] ds, \tag{25}$$

then

$$E\left(\sup_{T_2\leq s\leq T_2+t} |X_n(s)-X_i(s)|^2\right) \leq 4a_{i-1}+C_6\int_{T_2}^{T_2+t}\kappa\left[E\left(\sup_{T_2\leq s_1\leq s}|X_{n-1}(s_1-)-X_{i-1}(s_1-)|^2\right)\right]ds$$

$$\leq 4a_{i-1}+C_6\int_{T_2}^{T_2+t}\kappa(4C_1)ds_1$$

$$\leq 4a_{i-1}+C_5 t$$

This completes the proof of Lemma 6. □

Next, choose a positive number $0<\eta\leq T-T_2$ and a positive integer $p$ so that

$$C_6\kappa(4a_{p-1}+C_5 t)\leq C_5, \kappa_2(u)=C_6\kappa(u), t\in[0,\eta]. \tag{26}$$

Introduce the sequence of functions $\{\psi_k(t)\}_{k=1,2\cdots}$ and $\{\tilde{\psi}_{k,n}(t)\}_{k,n\geq 1}$, defined by

$$\psi_1(t)=4a_p+C_5 t,$$

$$\psi_{k+1}(t)=4a_{p+k}+\int_0^t \kappa_2(\psi_k(s))ds,$$

$$\tilde{\psi}_{k,n}(t)=E\left(\sup_{T_2\leq s\leq T_2+t}|X_{n+k}(s)-X_{p+k}(s)|^2\right)$$

**Lemma 7.** By Assumption 1,

$$\tilde{\psi}_{k,n}(t)\leq \psi_k(t)\leq \psi_{k-1}(t)\leq\cdots\leq \psi_1(t), t\in[0,\eta], \tag{27}$$

holds for all positive integer $k$.

**Proof:** We use Assumption 1, Eq. (25) and Lemma 6 to show it for $k=1,2$,

$$\tilde{\psi}_{1,n}(t)\leq 4a_p+C_6\int_{T_2}^{T_2+t}\kappa\left[E\left(\sup_{T_2\leq s_1\leq s}|X_n(s_1-)-X_p(s_1-)|^2\right)\right]ds$$

$$\leq 4a_p+\int_{T_2}^{T_2+t}\kappa_2(4a_{p-1}+C_5 s_1)ds_1=\psi_1(t),$$

and

$$\tilde{\psi}_{2,n}(t) \leq 4a_{p+1} + C_6 \int_{T_2}^{T_2+t} \kappa \left[ E\left( \sup_{T_2 \leq s_1 \leq s} |X_n(s_1-) - X_p(s_1-)|^2 \right) \right] ds$$

$$\leq 4a_{p+1} + \int_{T_2}^{T_2+t} \kappa_2(\tilde{\psi}_{1,n}(t)) ds_1 \leq 4a_{p+1} + \int_{T_2}^{T_2+t} \kappa_2(\psi_1(s_1)) ds_1$$

$$= \psi_2(t) \leq 4a_p + \int_{T_2}^{T_2+t} \kappa_2(4a_{p-1} + C_5 s_1) ds_1$$

$$\leq 4a_p + C_5 t = \psi_1(t), t \in [0, \eta].$$

Then we have proved

$$\tilde{\psi}_{2,n}(t) \leq \psi_2(t) \leq \psi_1(t).$$

Now assume that the assertion holds for $k \geq 2$. Then, by analogous argument, one obtains

$$\tilde{\psi}_{k+1,n}(t) \leq 4a_{p+k} + C_6 \int_{T_2}^{T_2+t} \kappa \left[ E\left( \sup_{T_2 \leq s_1 \leq s} |X_{n+k}(s_1-) - X_{p+k}(s_1-)|^2 \right) \right] ds$$

$$\leq 4a_{p+k} + \int_{T_2}^{T_2+t} \kappa_2(\tilde{\psi}_{k,n}(t)) ds_1$$

$$\leq 4a_{p+k} + \int_{T_2}^{T_2+t} \kappa_2(\psi_k(s_1)) ds_1 = \psi_{k+1}(t) \qquad (28)$$

$$\leq 4a_{p+k-1} + \int_{T_2}^{T_2+t} \kappa_2(\psi_{k-1}(s_1)) ds_1$$

$$= \psi_k(t), t \in [0, \eta].$$

Therefore, we obtain (27) for all $k$.

This completes the proof of Lemma 7. □

In terms of (27), we can define the function $\psi(t)$ by $\psi_k(t) \downarrow \psi(t), (k \to \infty)$ and yield that

$$\psi(0) = \lim_{k \to \infty} \psi_{k+1}(0) = \lim_{k \to \infty} a_{p+k} = 0.$$

Since $\psi(t)$ is a continuous function on $[0, \eta]$, by the definition of $\psi_{k+1}(t)$ and $\psi(t)$, we have

$$\psi(t) = \lim_{k \to \infty} \psi_{k+1}(t) = \lim_{k \to \infty} \left\{ 4a_{p+k} + \int_0^t \kappa_2(\psi_k(s)) ds \right\} = \int_0^t \kappa_2(\psi(s)) ds. \quad (29)$$

Then by $\psi(0) = 0$ and

$$\int_{0+} \frac{dq}{\kappa_2(q)} = \frac{1}{C_6} \int_{0+} \frac{dq}{\kappa(q)} = +\infty,$$

(29) implies $\psi(t) = 0, t \in [0, \eta]$.

Therefore, we obtain

$$\lim_{k \to \infty} \tilde{\psi}_{k,n}(t) = \lim_{k \to \infty} E\left( \sup_{0 \leq t \leq T_2 + \eta} |X_{n+k}(s) - X_{p+k}(s)|^2 \right)$$

$$\leq \lim_{k \to \infty} E\left( \sup_{0 \leq t \leq T_2} |X_{n+k}(s) - X_{p+k}(s)|^2 \right)$$

$$+ \lim_{k \to \infty} E\left( \sup_{T_2 \leq t \leq T_2 + \eta} |X_{n+k}(s) - X_{p+k}(s)|^2 \right)$$

$$\leq \lim_{k \to \infty} \psi_k(\eta) = \psi(\eta) = 0,$$

namely

$$\lim_{n,i \to \infty} E\left( \sup_{0 \leq t \leq T_2 + \eta} |X_n(t) - X_i(t)|^2 \right) = 0.$$

But, this conclusion is contradictory to the definition of $T_2$.

In other words, we have already shown that

$$\lim_{n,i \to \infty} E\left( \sup_{0 \leq t \leq T} |X_n(t) - X_i(t)|^2 \right) = 0. \tag{30}$$

The proof of the existence of solutions to Eq. (2) is complete. □

## 3.2 Uniqueness of solutions to SDEs

**Theorem 2.** Let $X(t)$ and $\tilde{X}(t)$ be two solutions to Eq. (2) on the same probability space such that $X(0) = \tilde{X}(0)$, then, under Assumption 1, the pathwise uniqueness holds for Eq. (2), $t \in [0, T]$.

**Proof:** By Cauchy-Schwarz inequality and Doob's martingale inequality, we observe that

$$E\left(\left|X(t)-\tilde{X}(t)\right|^2\right)$$

$$= E\left|\int_0^t \left[b(s, X(s)) - b(s, \tilde{X}(s))\right]ds + \int_0^t \left[\sigma(s, X(s)) - \sigma(s, \tilde{X}(s))\right]dB(s)\right.$$

$$\left. + \int_0^t \int_{|x|<c} \left[F(s, X(s-), x) - F(s, \tilde{X}(s-), x)\right]\tilde{N}(ds, dx)\right|^2$$

$$\leq 3TE\int_0^t \left|b(s, X(s-)) - b(s, \tilde{X}(s-))\right|^2 ds + 3\int_0^t \left\|\sigma(s, X(s-)) - \sigma(s, \tilde{X}(s-))\right\|^2 ds$$

$$+ 3\int_0^t \int_{|x|<c} \left|F(s, X(s-), x) - F(s, \tilde{X}(s-), x)\right|^2 v(dx)ds$$

$$\leq 3TE\int_0^t \left|b(s, X(s-)) - b(s, \tilde{X}(s-))\right|^2 + \left\|\sigma(s, X(s-)) - \sigma(s, \tilde{X}(s-))\right\|^2$$

$$+ \int_{|x|<c} \left|F(s, X(s-), x) - F(s, \tilde{X}(s-), x)\right|^2 v(dx)ds.$$

Noticing the Assumption 1, we have

$$E\left(\left|X(t)-\tilde{X}(t)\right|^2\right) \leq 3TE\int_0^t \lambda(s)\kappa\left(\left|X(s)-\tilde{X}(s)\right|^2\right)ds. \tag{31}$$

Since $\kappa(q)$ or $\kappa(q)^2/q$ is concave function, the above inequality (31) yields

$$E\left(\left|X(t)-\tilde{X}(t)\right|^2\right) \leq 3T\int_0^t \lambda(s)\kappa\left(E\left|X(s)-\tilde{X}(s)\right|^2\right)ds.$$

Then, noticing that $\lambda(t)$ an integrable function and $\int_{0+} \frac{du}{\kappa(q)} = +\infty$, the above inequality, as is well know, implies

$$E\left(\left|X(t)-\tilde{X}(t)\right|^2\right) = 0, t \in [0, T]. \tag{32}$$

Since $T$ is an arbitrary positive number, we obtain from this $X(t) \equiv \tilde{X}(t)$ for all $0 \leq t \leq T$.

The proof of the uniqueness of solutions of SDEs (2) is complete. □

## 4. Stability of solutions

**Definition 1.** A solution $X^\xi(t)$ of Eq.(2) with intial value $X(0) = \xi$ is said to be stable in mean square if for all $\varepsilon > 0$ there exists $\delta > 0$ such that when $E|\xi - \eta|^2 < \delta$,

$$E\left(\sup_{0\leq s\leq t}\left|X^{\xi}(s)-X^{\eta}(s)\right|^{2}\right)\leq\varepsilon,$$

where $X^{\eta}(t)$ is another solution of Eq.(2) with intial value $X(0)=\eta$.

In order to obtain the stability of solutions, we give two lemmas without details proof.

**Lemma 8.** (Bihari inequality) Let $T>0$ and $u_0 \geq 0$, $u(t), v(t)$ be continuous functions on $[0,T]$. Let $\kappa: R^+ \to R^+$ be a concave continuous and non-decreasing function such that $\kappa(q)>0$ for $q>0$. If

$$u(t)\leq u_0 + \int_0^t v(s)\kappa(u(s))ds, \text{ for all } t\in[0,T],$$

then

$$u(t)\leq G^{-1}\left(G(u_0)+\int_0^t v(s)ds\right),$$

for all $t\in[0,T]$ such that

$$G(u_0)+\int_0^t v(s)ds \in \text{Dom}(G^{-1}),$$

where $G(q)=\int_1^q \frac{1}{\kappa(s)}ds, q\geq 0$ and $G^{-1}$ is the inverse function of $G$.

**Lemma 9.** Let the assumptions of Lemma 8 hold and $v(t)\geq 0$ for $t\in[0,T]$. If for all $\varepsilon>0$, there exists $t_1\geq 0$ such that for $0\leq u_0 <\varepsilon$,

$$\int_{t_1}^{T}v(t)dt \leq \int_{u_0}^{\varepsilon}\frac{1}{\kappa(s)}ds,$$

holds. Then for every $t\in[t_1,T]$, the estimate $u(t)\leq\varepsilon$ holds.

Refer to [25, 26] to obtain proofs.

**Theorem 3.** Let $X^{\xi}(t)$ and $X^{\eta}(t)$ be two solutions of the Eq.(2) with intial value $\xi$ and $\eta$, respectively. The Assumption 1 is satisfied; then the solution of the Eq.(2) is said to be stable in mean square.

**Proof:**

$$E\left(\sup_{0\le s\le t}\left|X^{\xi}(s)-X^{\eta}(s)\right|^{2}\right)\le 4E|\xi-\eta|^{2}$$

$$+16tE\int_{0}^{t}\left|b\left(s,X^{\xi}(s-)\right)-b\left(s,X^{\eta}(s-)\right)\right|^{2}+\left\|\sigma\left(s,X^{\xi}(s-)\right)-\sigma\left(s,X^{\eta}(s-)\right)\right\|^{2}$$

$$+\int_{|x|<c}\left|F\left(s,X^{\xi}(s-),x\right)-F\left(s,X^{\eta}(s-),x\right)\right|^{2}v(dx)ds$$

$$\le 4E|\xi-\eta|^{2}+16t\left(\sup_{0\le s\le t}\lambda(s)\right)\int_{0}^{t}\kappa\left(E\sup_{0\le s_{1}\le s}\left|X^{\xi}(s_{1})-X^{\eta}(s_{1})\right|^{2}\right)ds.$$

Let

$$\kappa_{3}(q)=16T\left(\sup_{0\le s\le T}\lambda(s)\right)\kappa(q),$$

Thus, $\kappa_{3}(q)$ is obvious a monotone non-decreasing, continuous and concave function such that $\kappa_{3}(0)=0$ and $\int_{0+}\frac{1}{\kappa_{3}(q)}dq=\infty.$ So for any $\varepsilon>0$, $\varepsilon_{1}:=\frac{\varepsilon}{2}$, we have $\lim_{s\to 0}\int_{s}^{\varepsilon_{1}}\frac{1}{\kappa_{3}(q)}dq=\infty.$ Thus, there is a positive constant $\delta<\varepsilon_{1}$ such that

$$\int_{\delta}^{\varepsilon_{1}}\frac{1}{\kappa_{3}(q)}dq\ge T.$$

From the Lemma 9, let $u_{0}=4E|\xi-\eta|^{2}$, $v(t)=1$, $u(t)=E\left(\sup_{0\le s\le t}\left|X^{\xi}(s)-X^{\eta}(s)\right|^{2}\right).$

when $u_{0}\le\delta\le\varepsilon_{1}$, we have

$$\int_{u_{0}}^{\varepsilon_{1}}\frac{1}{\kappa_{3}(q)}dq\ge\int_{\delta}^{\varepsilon_{1}}\frac{1}{\kappa_{3}(q)}dq\ge T=\int_{0}^{T}v(t)dt.$$

So, for any $t\in[0,T]$, the estimate $u(t)\le\varepsilon$ holds.

This completes the proof. □

## Acknowledgments

This work was supported by the NSF of China (Grant Nos. 11372247, 11102157) and Shaanxi Province, Program for NCET, and SRF for ROCS, SEM, and NPU Foundation for New Faculties and Research Area Project and Garduate Strating Seed.